\documentclass[fractalfract,article,submit,moreauthors]{my-mdpi} 

\firstpage{1} 

\history{Received: 30-Dec-2020; Revised: 22-Feb-2021; Accepted: 04-March-2021}

\usepackage{amsmath,amssymb}
\usepackage{mathrsfs}

\allowdisplaybreaks 

\Title{Analysis of Hilfer Fractional Integro-Differential Equations with Almost Sectorial Operators}

\Author{Kulandhaivel Karthikeyan $^{1}$\orcidA{}, 
Amar Debbouche $^{2,3}$*\orcidB{} 
and Delfim F. M. Torres $^{3}$\orcidC{}}

\AuthorNames{Kulandhaivel Karthikeyan, Amar Debbouche and Delfim F. M. Torres}

\address{%
$^{1}$\quad  Department of Mathematics \& Centre for Research and Development, 
KPR Institute of Engineering and Technology, Coimbatore 641 407, Tamil Nadu, India;
karthi\_phd2010@yahoo.co.in\\
$^{2}$\quad Department of Mathematics, Guelma University, Guelma 24000, Algeria;
amar\_debbouche@yahoo.fr\\
$^{3}$\quad CIDMA---Center for Research and Development in Mathematics and Applications, 
Department of Mathematics, University of Aveiro, 3810-193 Aveiro, Portugal; delfim@ua.pt}

\corres{Correspondence: amar\_debbouche@yahoo.fr}

\abstract{We investigate a class of nonlocal integro-differential equations 
involving Hilfer fractional derivatives and almost sectorial operators. 
We prove existence results by applying Schauder's fixed point technique. 
Moreover, we show fundamental properties of the solution representation 
by discussing two cases related to the associated semigroup. For that, 
we consider compactness and noncompactness properties, respectively. 
Furthermore, an example is given to illustrate the obtained theory.}

\keyword{Hilfer fractional derivatives; mild solutions; almost sectorial operators; 
measure of non-compactness.}

\MSC{26A33; 47B12}

\begin{document}


\section{Introduction}

We consider nonlocal integro-differential equations involving Hilfer 
fractional derivatives and  almost sectorial operators:
\begin{equation}
\label{eq1}
D_{0^+}^{\alpha,\gamma}u(t)+\mathcal{A}u(t)
=g\left(t,u(t),\int^{t}_{0}k(t,s)f(s, u(s))ds\right),
\quad t\in(0,T]=\mathcal{J},
\end{equation}
\begin{equation}
\label{eq2}
I^{(1-\alpha)(1-\gamma)}_{0^+}[u(t)]\vert_{t=0}+h(u(t))=u_{0},
\end{equation}
where $D_{0^+}^{\alpha,\gamma}$ is the Hilfer fractional derivative of order $\alpha \in (0,1)$   
and type $\gamma \in [0,1]$. We assume that $\mathcal{A}$ is an almost sectorial operator 
on a Banach space $\mathcal{Y}$ with norm $\|\cdot\|$. Let $f: \mathcal{J} \times  \mathcal{Y}
\rightarrow  \mathcal{Y}$, $g: \mathcal{J} \times  \mathcal{Y}  \times  \mathcal{Y}  \rightarrow  \mathcal{Y}$ 
and $h: C(\mathcal{J}: \mathcal{Y})\rightarrow  \mathcal{Y}$ be given abstract functions, 
to be specified later. For brevity, we  take
\begin{align*}
\mathcal{B}u(t)=\int^{t}_{0}k(t,s)f(s,u(s))ds.
\end{align*}
   
During the last decades, mathematical modeling has been supported by the field of fractional calculus, 
with several successful results and fractional operators showing to be an excellent tool to describe 
the hereditary properties of various materials and processes. Recently, this combination has gained 
a lot of importance, mainly because fractional differential equations have become powerful tools 
in modeling several complex phenomena in numerous seemingly diverse and widespread fields of science 
and engineering, see, for instance, the basic text books \cite{book1,book2,book3,book4} and recent 
researches \cite{fr1,fr2,fr3}. In fact, abrupt changes, such as shocks, harvesting, 
or natural disasters, may happen in the dynamics of evolving processes. These short-term 
perturbations are often treated in the form of impulses. Recently, in so many published works, 
Hilfer fractional differential equations have received attention
\cite{h1,r2,f1,g1,a2,r3,a1,k1,t1,r1}.
   
In \cite{a1}, Jaiswal and Bahuguna study equations of Hilfer
fractional derivatives with almost sectorial operators in the abstract sense:
\begin{align*}
D^{\lambda,v}_{0^+}u(t)+\mathcal{A}u(t)=&g(t, u(t)),
\quad t\in (0,T],\\
\mathcal{I}^{(1-\lambda)(1-v)}_{0^+}u(0)=&u_{0}.
\end{align*}
We also  refer to the work in \cite{h1}, where Ahmed et al. study the question of existence  
for nonlinear Hilfer fractional differential equations with controls. Sufficient conditions 
are also established, where the time fractional derivative is the Hilfer derivative.
In \cite{z2}, Zhang and Zhou study fractional Cauchy problems 
with almost sectorial operators of the form
\begin{align*}
(^LD^{q}_{0+}x)(t)
=&Ax(t)+f(t, x(t)),~\text{ for almost all }~t\in [0,a],\\
(I^{(1-q)}_{0^+}x)(0)=&x_{0},
\end{align*}
where $^LD^{q}_{0+}$ is the Riemann--Liouville derivative of order $q$,  
$I^{(1-q)}_{0^+}$ is the Riemann--Liouville integral of order $1-q$, $0<q<1$, 
$A$ is an almost sectorial operator on a complex Banach space, and $f$ is a given function. 
Motivated by these results, here we extend the previous available results 
of the literature to a class of Hilfer fractional integro-differential equations 
in which the closed operator is almost sectorial. Moreover, we also  
consider both compactness and noncompactness cases of the semigroup operator. 

The paper is structured as follows. In Section~\ref{sec:2}, we present necessary information  
about the Hilfer derivative, almost sectorial operators, measure of non-compactness,  
mild solutions of equations \eqref{eq1}--\eqref{eq2}, along with some useful definitions, 
results, and lemmas. We discuss fundamental results for mild solutions for the equations
\eqref{eq1}--\eqref{eq2} in Section~\ref{sec:3}. In Section~\ref{sec:4}, we prove the solvability
question in two cases, when associated semigroup is compact and noncompact, respectively. 
An example is then given in Section~\ref{sec:5}, to illustrate our main results.
We end with Section~\ref{sec:6} of conclusions.        


\section{Preliminaries}
\label{sec:2}

In this section we recall some necessary theory that will be used throughout 
the work in order to get the new results. 


\subsection{Fractional derivatives}

We start by a short introduction of the main definitions in fractional calculus \cite{r1,z1}. 

\begin{Definition}
\label{Definition2.1}
The left-sided Riemann--Liouville fractional integral of order $\alpha>0$ with lower limit $a$ 
for a function $h: [a, +\infty)\to \mathbb{R}$ is defined as
$$
I^{\alpha}_{a^{+}}h(t)=\frac{1}{\Gamma(\alpha)}\int_{a}^{t}(t-s)^{\alpha-1}h(s)ds,
$$
provided the right hand side is defined a.e. on $[a, +\infty)$. 
\end{Definition}

\begin{Remark}
If $a=0$, then we write
$I^{\alpha}_{0^{+}}f(t)=(g_{\alpha}*f)(t)$, where
$$
g_{\alpha}(t):=\left\{
\begin{array}{ll}
\frac{1}{\Gamma(\alpha)}t^{\alpha-1},& \mbox{$t>0$},\\
0, & \mbox{$t\leq 0$},
\end{array}\right.
$$
and, as usual, $*$ denotes the convolution of functions. 
Note that $\lim\limits_{\alpha\rightarrow 0^+}g_{\alpha}(t)=\delta(t)$ 
with $\delta$ the delta Dirac function.
\end{Remark}

\begin{Definition}
\label{Definition2.2}
The left-sided Riemann--Liouville fractional derivative of order $\alpha>0$, 
$n-1\leq\alpha<n$, $n\in \mathbb{N}$, for a function 
$h: [a, +\infty)\to \mathbb{R}$, is defined by
$$
^{L}D^{\alpha}_{a^{+}}h(t)=\frac{1}{\Gamma(n-\alpha)}
\frac{d^{n}}{dt^{n}}\int_{a}^{t}\frac{h(s)}{(t-s)^{\alpha+1-n}}ds,
\quad t>a,
$$
provided the right hand side is defined a.e. on $[a, +\infty)$.
\end{Definition}

\begin{Definition}
\label{Definition2.3}
The left-sided Caputo's fractional derivative of order $\alpha>0$, 
$n-1<\alpha<n$, $n\in \mathbb{N}$, for a function $h: [a, +\infty)\to \mathbb{R}$, 
is defined by
$$
^{C}D^{\alpha}_{a^{+}}h(t)=\frac{1}{\Gamma(n-\alpha)}
\int_{a}^{t}\frac{h^{(n)}(s)}{(t-s)^{\alpha+1-n}}ds=I^{n-\alpha}_{a^{+}}h^{(n)}(t),
\quad t>a,
$$
provided the right hand side is defined a.e. on $[a, +\infty)$.
\end{Definition}

\begin{Definition}
\label{Definition2.4}
The left-sided Hilfer fractional derivative of order $0<\alpha<1$ 
and type $\gamma\in[0, 1]$, of a function $h: [a, +\infty)\to \mathbb{R}$, 
is defined as
$$
D^{\alpha,\gamma}_{a^{+}}h(t)
=\left[I^{(1-\alpha)\gamma}_{a^{+}}D\left(I^{(1-\alpha)(1-\gamma)}_{a^{+}}h\right)\right](t).
$$
\end{Definition}

\begin{Remark}
\label{Remark2.1}
(i) If $\gamma=0$, $0<\alpha<1$, and $a=0$, 
then the Hilfer fractional derivative corresponds 
to the classical Riemann--Liouville fractional derivative:
$$
D^{\alpha,0}_{0^{+}}h(t)=\frac{d}{dt}I^{1-\alpha}_{0^{+}}h(t)=~^{L}D^{\alpha}_{0^{+}}h(t).
$$
(ii) If $\gamma=1$, $0<\alpha<1$, and $a=0$, then the Hilfer fractional derivative 
corresponds to the classical Caputo fractional derivative:
$$
D^{\alpha,1}_{0^{+}}h(t)=I^{1-\alpha}_{0^{+}}\frac{d}{dt}h(t)=~^{C}D^{\alpha}_{0^{+}}h(t).
$$
\end{Remark}


\subsection{Measure of non-compactness}

The motivation to consider our problem can be found in \cite{z2,z1}.
Here we generalize the results in \cite{z2,z1}.
Let $\mathcal{L}$ $\subset$  $\mathcal{Y}$ be bounded. 
The Hausdorff measure of non-compactness  is considered as 
\begin{align*}
\Theta(\mathcal{L})
=\inf\left\{  \theta >0 \quad \text{such that} \quad \mathcal{L}
\subset \bigcup_{j=1}^{m} B_{ \theta}(x_{j}),\quad 
\text{where}\quad x_{j} \in  \mathcal{Y}, \ m \in N \right\} 
\end{align*}
while the Kurtawoski measure of noncompactness $\Phi$ on a bounded 
set $\mathcal{B}\subset  \mathcal{Y} $ is given by
\begin{align*}
\Phi(\mathcal{L})=\inf\left\{ \epsilon >0 \quad \text{such that } 
\quad \mathcal{L}\subset \bigcup_{j=1}^{m}M_{j} 
\quad \text{and}\quad \text{diam}(M_{j})\leq \epsilon \right\}
\end{align*}
with the following properties:
\begin{enumerate}
\item $\mathcal{L}_{1} \subset \mathcal{L}_{2}$ 
gives $\Theta(\mathcal{L}_{1} )\leq \Theta(\mathcal{L}_{2} )$, 
where $\mathcal{L}_{1},\mathcal{L}_{2}$ 
are bounded subsets of $\mathcal{Y}$;

\item $\Theta(\mathcal{L})=0$ if and only if 
$\mathcal{L}$ is relatively compact in  $\mathcal{Y}$;

\item $\Theta(\{z\}\bigcup \mathcal{L})=\Theta(\mathcal{L}) $ for all 
$z\in  \mathcal{Y}$  $\mathcal{L}\subseteq  \mathcal{Y}$;

\item $\Theta(\mathcal{L}_{1}\bigcup \mathcal{L}_{2} ) 
\leq \max\{ \Theta(\mathcal{L}_{1}), \Theta(\mathcal{L}_{2})\}$;

\item $\Theta(\mathcal{L}_{1}+\mathcal{L}_{2} ) 
\leq  \Theta(\mathcal{L}_{1})+\Theta(\mathcal{L}_{2})$;

\item $\Theta(r\mathcal{L}) \leq |r| \Theta(\mathcal{L})$ for $r \in \mathbb{R}$.
\end{enumerate}		
Let $\mathcal{M} \subset C(I, \mathcal{Y} )$ and  
$\mathcal{M}(r)=\left\{ \upsilon(r) 
\in  \mathcal{Y} | \upsilon \in \mathcal{M} \right\}$. One defines
\begin{align*}
\int_{0}^{t}\mathcal{M}(r)dr
:=\left\{ \int_{0}^{t}\upsilon(r)dr|\upsilon\in \mathcal{M} \right\},
\quad t \in \mathcal{J}.
\end{align*}
	
\begin{Proposition}
\label{Proposition2.2}  
If $\mathcal{M}\subset C(\mathcal{J},  \mathcal{Y} )$ is equicontinuous and bounded, 
then $t\rightarrow \Theta(\mathcal{M}(t))$ is continuous on $I$. Also,
\begin{align*}
\Theta(\mathcal{M})=\max\left\{\Theta(\mathcal{M}(t)),
\Theta\left(\int_{0}^{t}\upsilon(r)dr\right)\right\}
\leq \int_{0}^{t}\Theta(\upsilon(r))dr,\quad \text{for}\quad t\in I. 
\end{align*} 
\end{Proposition}

\begin{Proposition}
\label{Proposition2.3} 
Let $\{\upsilon_{n}:\mathcal{J}\rightarrow  \mathcal{Y} , n \in \mathbb{N}\}$ 
be Bochner integrable functions. This implies that 
$\|\upsilon_{n}\|\leq m(t)$ a.e. for $n\in \mathbb{N}$ and
$m\in L^{1}(I,R^{+})$. Then, 
$\xi(t)=\Theta\left(\{\upsilon_{n}(t)\}_{n=1}^{\infty}\right) 
\in L^{1}(I,R^{+})$ and satisfies 
\begin{align*}
\Theta\left(\left\{\int_{0}^{t}\upsilon_{n}(r)dr
: n\in \mathbb{N}\right\}\right)
\leq 2 \int_{0}^{t}\xi(r)dr.
\end{align*}
\end{Proposition}

\begin{Proposition}
\label{Proposition2.4} 
Let $\mathcal{M}$ be a bounded set. Then, for any $ \theta>0$, 
there exists a sequence $\{\upsilon_{n}\}_{n=1}^{\infty}\subset \mathcal{M}$ such that
\begin{align*}
\Theta(\mathcal{M}) \leq 2\Theta \{\upsilon_{n}\}_{n=1}^{\infty}\ +  \theta.
\end{align*}
\end{Proposition}


\subsection{Almost sectorial operators}

Let $0<\beta<\pi $ and $-1<\beta<0$. We define 
$S_{\beta}^{0}:=\{\upsilon\in C\setminus \{0\}~ \text{ such that}~|\arg \upsilon|<\beta\}$ 
and its closure by $S_{\beta}$, such that $ S_{\beta}=\{\upsilon\in C\setminus \{0 \}$ 
with $|\arg \upsilon|\leq\beta\} \bigcup \{0\}$.

\begin{Definition}[See \cite{p1}]
\label{Definition2.6}  
For $-1<\beta<0$ and $0<\omega<\frac{\pi}{2}$, we define 
$\left\{ \Theta_{\omega}^{\beta}\right\}$ as the family 
of all closed and linear operators 
$\mathcal{A}:D(\mathcal{A})\subset 
\mathcal{Y} \rightarrow  \mathcal{Y}$ such that
\begin{enumerate}
\item $\sigma (\mathcal{A})$ is contained in  $S_{\omega}$;

\item for all $ \beta \in (\omega,\pi)$ there exists $M_{\beta}$ such that
\begin{align*}
\|\mathcal{R}(z,\mathcal{A})\|_{L(\mathcal{Y})}\leq M_{\beta}|z|^{\beta},
\end{align*}
where $\mathcal{R}(z,\mathcal{A})=(zI-\mathcal{A})^{-1}$ is the resolvent operator 
and $\mathcal{A} \in \Theta_{\omega}^{\beta}$ is said to be an almost 
sectorial operator on $\mathcal{Y}$.
\end{enumerate}
\end{Definition} 

\begin{Proposition}[See \cite{p1}]
\label{Proposition2.5} 
Let $\mathcal{A} \in \Theta_{\omega}^{\beta}$ for $-1<\beta<0$ 
and $0<\omega<\frac{\pi}{2}$. Then the following properties hold:
\begin{enumerate}
\item $ Q(t)$ is analytic and 
$\frac{d^{n}}{dt^{n}} Q(t)=(-\mathcal{A}^{n} Q(t)(t\in S_{\frac{\pi}{2}}^{0})$;

\item $Q(t+s)= Q(t)  Q(s)$ $\forall$ $t,s \in S_{\frac{\pi}{2}}^{0}$;

\item $ \| Q(t)\|_{L( \mathcal{Y} )}\leq C_{0}t^{-\beta-1}(t>0)$, 
where $C_{0}=C_{0}(\beta)>0$ is a constant;

\item if $\sum_{  Q}=\{x \in  \mathcal{Y} ~~:lim_{t\rightarrow 0_{+}}  Q(t)x=x \}$,
then $D(\mathcal{A}^{\theta})\subset \sum_{  Q}$ if $\theta > 1+\beta$;
 
\item $\mathcal{R}(r,-\mathcal{A})
=\int_{0}^{\infty}e^{-rs}  Q(s)ds$, $r\in \mathbb{C}$ with $Re(r)>0$.
\end{enumerate}
\end{Proposition}

We use the following Wright-type function \cite{z1}:
\begin{align*}
M_{\alpha}(\theta)=\sum_{n\in \mathbb{N}}\frac{(-\theta)^{n-1}}{\Gamma (1-{\alpha}n)(n-1)!},
\quad \theta \in \mathbb{C}.
\end{align*}
For $-1<\sigma<\infty,~ r>0$, the following properties hold: 
\begin{enumerate}
\item[(A1)] $M_{\alpha}(\theta)\geq 0$, $t>0$;
		
\item[(A2)] $\int_{0}^{\infty}{\theta}^{\sigma}
M_{\alpha}d\theta=\frac{\Gamma(1+\sigma)}{\Gamma(1+\alpha \sigma)}$;
		
\item[(A3)] $\int_{0}^{\infty} \frac{\alpha}{\theta^{\alpha+1}}
e^{-r\theta}M_{\alpha}(\frac{1}{\theta^{\alpha}})d\theta=e^{-r^\alpha}$.
\end{enumerate}
The characteristic operators  
$\{ S_{\alpha}(t)\}|_{t\in S_{\frac{\pi}{2}-w}^{0}}$ 
and $\{ T_{\alpha}(t)\}|_{t\in S_{\frac{\pi}{2}-w}^{0}}$ are defined by		
\begin{align*}
S_{\alpha}(t):=\int_{0}^{\infty}M_{\alpha}( \theta) Q(t^{\alpha} \theta)d \theta
\end{align*}
and
\begin{align*}
T_{\alpha}(t):=\int_{0}^{\infty}\alpha  
\theta M_{\alpha}( \theta) Q(t^{\alpha} \theta)d \theta.
\end{align*}

\begin{Theorem}[See Theorem 4.6.1 of \cite{z1}]
\label{Theorem2.1} 
For each fixed	$t \in S_{\frac{\pi}{2}-\omega}^{0}$, $S_{\alpha}(t)$ 
and $ T_{\alpha}(t)$ are bounded linear operators on $\mathcal{Y}$. 
Moreover,
\begin{align*}
\| S_{\alpha}(t)\|\leq  C_1t^{-\alpha(1+\beta)},
\quad \| T_{\alpha}(t)\|\leq  C_2t^{-\alpha(1+\beta)},
\quad t>0,
\end{align*}
where $ C_1$ and $ C_2$ are constants dependent on $\alpha$ and $\beta$.
\end{Theorem}				

\begin{Theorem}[See \cite{z1}]
\label{Theorem2.2}   
The operators $S_{\alpha}(t)$ and $ T_{\alpha}(t)$ are continuous 
in the uniform operator topology, for $t>0$. For $s>0$, 
the continuity is uniform on $[s,\infty]$.
\end{Theorem}
			
Define $\Omega_{r}(\mathcal{J}):=
\lbrace y\in C(\mathcal{J}, \mathcal{Y})\vert~ \Vert y\Vert\leq r\rbrace$. 
Our main results are proved under the following hypotheses:
\begin{enumerate}

\item[(H1)] For  $t\in \mathcal{J}$, 
$g(t,\cdot,\cdot): \mathcal{Y}\times  \mathcal{Y} \rightarrow  \mathcal{Y}$ 
and $f(t,\cdot): \mathcal{Y} \rightarrow  \mathcal{Y}$ are continuous functions and, 
for each $u\in  C(\mathcal{J}, \mathcal{Y} )$, 
$g(\cdot,u,\mathcal{B}u): \mathcal{J}\rightarrow  \mathcal{Y} $ 
and $f(\cdot, u): \mathcal{J}\rightarrow  \mathcal{Y} $ are strongly measurable.
 
\item[(H2)] There exist functions $k_1, k_2\in L^{1}(\mathcal{J},\mathbb{R}^{+})$ satisfying 
$\Vert g(t,\cdot,\cdot)\Vert\leq k_1(t)+k_2(t)e^{-\delta t}$ for all $u\in \Omega_{r}(\mathcal{J})$ 
and almost all $t$ on $\mathcal{J}$ and  
$$
I_{0^+}^{-\alpha\beta}[k_1(t)+k_2(t)e^{-\delta t}]
\in  C(\mathcal{J},\mathbb{R}),~\lim\limits_{t\to 0^+}~t^{(1+\alpha\beta)(1-\gamma)}
I_{0^+}^{-\alpha\beta}[k_1(t)+k_2(t)e^{-\delta t}]=0.
$$
 	
\item[(H3)] Function $h: C(\mathcal{J}, \mathcal{Y})\rightarrow  \mathcal{Y}$ 
is completely continuous and there exists a positive constant $k$ such that $\Vert h(u)\Vert\leq k$.

\item[(H4)] We assume that
\begin{align*}
\sup_{[0,T]}(t^{(1+\alpha\beta)(1-\gamma)}\| S_{\alpha,\gamma}(t)[u_{0}+k]\|
+t^{(1+\alpha\beta)(1-\gamma)}\int_{0}^{t}(t-r)^{-\alpha\beta-1}[k_1(r)+k_2(r)e^{-\delta r}]dr)\leq r,
\end{align*}
for $r>0$, $u_{0}\in D( \mathcal{A}^{\theta})$, and $\theta>1+\beta$, where 
$S_{\alpha,\gamma}(t)=I_{0^+}^{\gamma(1-\alpha)}t^{\alpha -1} T_{\alpha}(t)$.
\end{enumerate}

For the next two lemmas, we refer to \cite{r2,t1}. 
 
\begin{Lemma}[See \cite{r2,t1}]
\label{Lemma2.1} 
The fractional Cauchy problem \eqref{eq1}--\eqref{eq2} 
is equivalent to the integral equation
\begin{align}
\label{eq3}
u(t)=&\frac{[u_{0}-h(u(t))]}{\Gamma(\gamma(1-\alpha)+\alpha)}t^{(1-\alpha)(\gamma-1)}
+\frac{1}{\Gamma(\alpha)}\int_{0}^{1}(t-r)^{\alpha-1}[- \mathcal{A} u(r)+g(r,u(r),(\mathcal{B}u)r)]dr,
\quad t\in \mathcal{J}.
\end{align}
\end{Lemma}

\begin{Lemma}[See \cite{r2,t1}]
\label{Lemma2.2} 
If $u$ is a solution to the integral equation \eqref{eq3}, then it satisfies
\begin{align*}
u(t)= S_{\alpha,\gamma}(t)[u_{0}-h(u(t))]+\int_{0}^{t} R_{\alpha}(t-r)g(r,u(r),(\mathcal{B}u)r)dr,
\end{align*}
where $S_{\alpha,\gamma}(t)=I_{0^+}^{\gamma(1-\alpha)} R_{\alpha}(t)$ 
with $R_{\alpha}=t^{\alpha-1} T_{\alpha}(t)$.
\end{Lemma}

\begin{Definition}
\label{Definition2.7} 
By a mild solution of the Cauchy problem \eqref{eq1}--\eqref{eq2}, 
we mean a function $u\in C(\mathcal{J}, \mathcal{Y})$ that satisfies
\begin{align*}
u(t)= S_{\alpha,\gamma}(t)[u_{0}-h(u(t))]
+\int_{0}^{t} R_{\alpha}(t-r)g(r,u(r),(\mathcal{B}u)r)dr,
\quad t\in \mathcal{J}. 
\end{align*}
\end{Definition}

We define operator  
$\mathcal{P}:\Omega_{r}(\mathcal{J})\rightarrow \Omega_{r}(\mathcal{J})$ as
\begin{align*}
(\mathcal{P}u)(t):= S_{\alpha,\gamma}(t)[u_{0}-h(u(t))]
+\int_{0}^{t}(t-r)^{\alpha-1} T_{\alpha}(t-r)g(r,u(r),(\mathcal{B}u)r)dr.
\end{align*}

\begin{Lemma}[See \cite{a1}]
\label{Lemma2.3} 
The operators $ R_{\alpha}(t)$ and $ S_{\alpha,\gamma}(t)$ are bounded linear operators 
on $\mathcal{Y}$ for every fixed $t\in S_{\frac{\pi}{2}-\omega}^{0}$. 
Also, for $t>0$, we have 	
\begin{align*}
\| R_{\alpha}(t)x\|\leq  C_2t^{-1-\alpha\beta}\|x\|,
\quad \| S_{\alpha,\gamma}(t)x\|
\leq \frac{\Gamma(-\alpha\beta)}{\Gamma(\gamma(1-\alpha)-\alpha\beta)}
C_2t^{\gamma(1-\alpha)-\alpha\beta-1}\|x\|.
\end{align*}
\end{Lemma}

\begin{Proposition}[See \cite{a1}]
\label{Proposition2.6}  
The operators $ R_{\alpha}(t)$ and $ S_{\alpha,\gamma}(t)$ 
are strongly continuous for $t>0$. 
\end{Proposition}


\section{Auxiliary results}
\label{sec:3}

Now, we are in position to start our original contributions.

\begin{Theorem}
\label{Theorem3.1}
Let $\mathcal{A} \in \Theta_{\omega}^{\beta}$ for 
$-1<\beta<0$ and $0<\omega<\frac{\pi}{2}$. Assuming that $(H1)$--$(H4)$ are satisfied, 
then the operator $\{\mathcal{P}y~:~y\in \Omega_{r}(\mathcal{J})\}$ is equicontinuous, 
provided $u_{0}\in  \mathcal{D}(\mathcal{A}^{\theta})$ with $\theta>1+\beta$.
\end{Theorem}

\begin{proof}
For $y\in \Omega_{r}(\mathcal{J})$ and $t_{1}=0<t_{2}\leq T$, we have
\begin{align*}
\Big\| \mathcal{P}y(t_{2})&- \mathcal{P}y(0)\Big\|\\
&=\Big\|t_{2}^{(1+\alpha\beta)(1-\gamma)}\Big(
S_{\alpha,\gamma}(t_{2})[u_{0}-h(u(t))]+\int_{0}^{t_{2}}(t_{2}-r)^{\alpha-1}
T_{\alpha}(t_{2}-r)g(r,u(r),(\mathcal{B}u)r)dr\Big)\Big\|\\
&\leq \Big\|t_{2}^{(1+\alpha\beta)(1-\gamma)}
S_{\alpha,\gamma}(t_{2})\Big\|(u_{0}+k)\\
&\quad +\Big\|t_{2}^{(1+\alpha\beta)(1-\gamma)}
\int_{0}^{t_{2}}(t_{2}-r)^{\alpha-1} 
T_{\alpha}(t_{2}-r)g(r,u(r),(\mathcal{B}u)r)dr\Big\|
\rightarrow 0, \quad \text{ as } t_{2}\rightarrow 0.
\end{align*}
Now, let $0<t_{1}<t_{2}\leq T$. One has
\begin{align*}
\Big\| \mathcal{P}y(t_{2})- \mathcal{P}y(t_{1})\Big\|
\leq& \Big\|t_{2}^{(1+\alpha\beta)(1-\gamma)} 
S_{\alpha,\gamma}(t_{2})[u_{0}-h(u(t))]-t_{1}^{(1+\alpha\beta){1-\gamma}}
S_{\alpha,\gamma}(t_{1})[u_{0}-h(u(t))]\Big\|\\
&+\Big\|t_{2}^{(1+\alpha\beta)(1-\gamma)}
\int_{0}^{t_{2}}(t_{2}-r)^{\alpha-1} T_{\alpha}(t_{2}-r)g(r,u(r),(\mathcal{B}u)r)dr\\
&-t_{1}^{(1+\alpha\beta)(1-\gamma)}\int_{0}^{t_{1}}(t_{1}-r)^{\alpha-1}
T_{\alpha}(t_{1}-r)g(r,u(r),(\mathcal{B}u)r)dr\Big\|.
\end{align*}
Using the triangle inequality, we get
\begin{align*}
\Big\| \mathcal{P}y(t_{2})- \mathcal{P}y(t_{1})\Big\|
\leq & \Big\|t_{2}^{(1+\alpha\beta)(1-\gamma)}
S_{\alpha,\gamma}(t_{2})[u_{0}-h(u(t))]-t_{1}^{(1+\alpha\beta)(1-\gamma)}
S_{\alpha,\gamma}(t_{1})[u_{0}-h(u(t))]\Big\|\\
&+\Big\|t_{2}^{(1+\alpha\beta)(1-\gamma)}
\int_{t_{1}}^{t_{2}}(t_{2}-r)^{\alpha-1} T_{\alpha}(t_{2}-r)g(r,u(r),(\mathcal{B}u)r)dr\Big\|\\
&+\Big\|t_{2}^{(1+\alpha\beta)(1-\gamma)}\int_{0}^{t_{1}}(t_{2}-r)^{\alpha-1}
T_{\alpha}(t_{2}-r)g(r,u(r),(Bu)r)dr\\
&-t_{1}^{(1+\alpha\beta)(1-\gamma)}\int_{0}^{t_{1}}(t_{1}-r)^{\alpha-1}
T_{\alpha}(t_{2}-r)g(r,u(r),(\mathcal{B}u)r)dr\Big\|\\&
+\Big\|t_{1}^{(1+\alpha\beta)(1-\gamma)}\int_{0}^{t_{1}}(t_{1}-r)^{\alpha-1}
T_{\alpha}(t_{2}-r)g(r,u(r),(\mathcal{B}u)r)dr\\
&-t_{1}^{(1+\alpha\beta)(1-\gamma)}\int_{0}^{t_{1}}(t_{1}-r)^{\alpha-1}
T_{\alpha}(t_{1}-r)g(r,u(r),(\mathcal{B}u)r)dr\Big\|\\
=:& \mathfrak{I}_{1}+\mathfrak{I}_{2}+\mathfrak{I}_{3}+\mathfrak{I}_{4}.
\end{align*}
From the strong continuity of $ S_{\alpha,\gamma}(t)$, 
we have $\mathfrak{I}_{1}\rightarrow 0$ as $t_{2}\rightarrow t_{1}$. Also,
\begin{align*}
\mathfrak{I}_{2}\leq & 
C_2t_{2}^{(1+\alpha\beta)(1-\gamma)}
\int_{t_{1}}^{t_{2}}(t_{2}-r)^{-\alpha\beta-1}[k_1(r)+k_2(r)e^{-\delta r}]dr\\ 
\leq & C_2\Big|t_{2}^{(1+\alpha\beta)(1-\gamma)}
\int_{0}^{t_{2}}(t_{2}-r)^{-\alpha\beta-1}[k_1(r)+k_2(r)e^{-\delta r}]dr\\
&\quad -t_{2}^{(1+\alpha\beta)(1-\gamma)}
\int_{0}^{t_{1}}(t_{1}-r)^{-\alpha\beta-1}[k_1(r)+k_2(r)e^{-\delta r}]dr \Big|\\ 
\leq &  C_2\int_{0}^{t_{1}} \Big|t_{1}^{(1+\alpha\beta)(1-\gamma)}(t_{1}-r)^{-\alpha\beta-1}
-t_{2}^{(1+\alpha\beta)(1-\gamma)}(t_{2}-r)^{-\alpha\beta-1}\Big|[k_1(r)+k_2(r)e^{-\delta r}]dr.
\end{align*}
Then, by using $(H2)$ and the dominated convergence theorem,
$\mathfrak{I}_{2}\rightarrow 0$ as $ t_{2}\rightarrow t_{1}$. Since
\begin{align*}
\mathfrak{I}_{3}\leq 
C_2\int_{0}^{t_{1}}(t_{2}-r)^{-\alpha-\alpha\beta}\Big|
t_{2}^{(1+\alpha\beta)(1-\gamma)}(t_{2}-r)^{\alpha-1}
-t_{1}^{(1+\alpha\beta)(1-\gamma)}(t_{1}-r)^{\alpha-1}
\Big|[k_1(r)+k_2(r)e^{-\delta r}]dr,
\end{align*} 
\begin{align*}
(t_{2}-r)&^{-\alpha-\alpha\beta}\Big|t_{2}^{(1+\alpha\beta)(1-\gamma)}(t_{2}
-r)^{\alpha-1}-t_{1}^{(1+\alpha\beta)(1-\gamma)}(t_{1}-r)^{\alpha
-1}\Big|[k_1(r)+k_2(r)e^{-\delta r}]\\
\leq & t_{2}^{(1+\alpha\beta)(1-\gamma)}(t_{2}-r)^{\alpha-1}[k_1(r)+k_2(r)e^{
-\delta r}]+t_{1}^{(1+\alpha\beta)(1-\gamma)}(t_{1}-r)^{\alpha
-1}[k_1(r)+k_2(r)e^{-\delta r}]\\
\leq & 2t_{1}^{(1+\alpha\beta)(1-\gamma)}(
t_{1}-r)^{\alpha-1}[k_1(r)+k_2(r)e^{-\delta r}]
\end{align*}
and $\int_{0}^{t_{1}}2t_{1}^{(1+\alpha\beta)(1-\gamma)}(t_{1}
-r)^{\alpha-1}[k_1(r)+k_2(r)e^{-\delta r}]dr$ exists, 
we obtain $I_{3}\rightarrow 0$ as $t_{2}\rightarrow t_{1}$.
For $\epsilon >0$, we have
\begin{align*}
\mathfrak{I}_{4}
=&\Big\|\int_{0}^{t_{1}}t_{1}^{(1+\alpha\beta)(1-\gamma)}[ T_{\alpha}(t_{2}-r)
- T_{\alpha}(t_{1}-r)](t_{1}-r)^{\alpha-1}g(r,u(r),(\mathcal{B}u)r)dr\Big\|\\
\leq & \int_{0}^{t_{1}-\epsilon}t_{1}^{(1+\alpha\beta)(1-\gamma)}\Big\| T_{\alpha}(t_{2}-r)
- T_{\alpha}(t_{1}-r)\Big\|_{L(\mathcal{Y})}(t_{1}-r)^{\alpha-1}[k_1(r)+k_2(r)e^{-\delta r}]\\
&+ \int_{t_{1}}^{t_{1}}t_{1}^{(1+\alpha\beta)(1-\gamma)}\Big\| T_{\alpha}(t_{2}-r)
- T_{\alpha}(t_{1}-r)\Big\|_{L(\mathcal{Y})}(t_{1}-r)^{\alpha-1}[k_1(r)+k_2(r)e^{-\delta r}]\\ 
\leq & t_{1}^{(1+\alpha\beta)(1-\gamma)}\int_{0}^{t_{1}}(t_{1}-r)^{\alpha-1}[k_1(r)
+k_2(r)e^{-\delta r}]dr~~\sup_{s\in [0,t_{1}-\epsilon]}~~\Big\| T_{\alpha}(t_{2}-r)
- T_{\alpha}(t_{1}-r)\Big\|_{L(\mathcal{Y})}\\ 
&+ C_2\int_{t_{1}}^{t_{1}}t_{1}^{(1+\alpha\beta)(1-\gamma)}((
t_{2}-r)^{-\alpha-\alpha\beta}+(t_{1}-r)^{-\alpha-\alpha\beta})(
t_{1}-r)^{\alpha-1}[k_1(r)+k_2(r)e^{-\delta r}]dr\\
\leq & t_{1}^{(1+\alpha\beta)(1-\gamma)+\alpha(1+\beta)}
\int_{0}^{t_{1}}(t_{1}-r)^{-\alpha\beta-1}[k_1(r)
+k_2(r)e^{-\delta r}]dr~\sup_{s
\in [0,t_{1}-\epsilon]}\Big\| T_{\alpha}(t_{2}-r)
- T_{\alpha}(t_{1}-r)\Big\|_{L(\mathcal{Y})}\\ 
& +2 C_2\int_{t_{1}-\epsilon}^{t_{1}}t_{1}^{(1+\alpha\beta)(1-\gamma)}(t_{1}
-r)^{-\alpha\beta-1}[k_1(r)+k_2(r)e^{-\delta r}]dr.
\end{align*}
Since $ T_{\alpha}(t)$ is uniformly continuous 
and $\lim_{t_{2}\rightarrow t_{1}}\mathfrak{I}_{2}=0$, 
then $\mathfrak{I}_{4}\rightarrow 0$ as $t_{2}\rightarrow 
t_{1}$, independent of $y\in \Omega_{r}(\mathcal{J})$.
Hence, $\Big\| \mathcal{P}y(t_{2})- \mathcal{P}y(t_{1})\Big\|\rightarrow 0$, 
independently of $y\in \Omega_{r}(\mathcal{J})$  as $t_{2}\rightarrow t_{1}$. 
Therefore, $\{\mathcal{P}y~:~y\in \Omega_{r}(\mathcal{J})\}$ is equicontinuous.
\end{proof}
 
\begin{Theorem}
\label{Theorem3.2}
Let $-1<\beta<0$, $0<\omega<\frac{\pi}{2}$, 
and $\mathcal{A} \in \Theta_{\omega}^{\beta}$. 
Then, under hypotheses $(H1)$--$(H4)$, the operator 
$\{\mathcal{P}y~:~y\in \Omega_{r}(\mathcal{J})\}$ is continuous 
and bounded, provided 
$u_{0}\in  \mathcal{D}(\mathcal{A}^{\theta})$ with $\theta>1+\beta $. 
\end{Theorem}

\begin{proof}
We verify that $\mathcal{P}$ maps $\Omega_{r}(\mathcal{J})$ into itself. 
Taking $y\in \Omega_{r}(\mathcal{J}) $ and defining 
$$
u(t):=t^{-(1+\alpha\beta)(1-\gamma)}y(t),
$$ 
we have $u\in \Omega_{r}(\mathcal{J})$. Let $ t \in [0,T]$,
\begin{align*}
\|\mathcal{P}\|\leq \|t^{(1+\alpha\beta)(1-\gamma)}
S_{\alpha,\gamma}(t)[u_{0}-h(u(t))]\|+t^{(1+\alpha\beta)(1-\gamma)}
\Big\|\int_{0}^{t}(t-r)^{\alpha-1} 
T_{\alpha}(t-r)g(r,u(r),(\mathcal{B}u)r)dr\Big\|.
\end{align*}
From $(H2)$--$(H4)$, we get
\begin{align*}
\| \mathcal{P}y(t)\|
&\leq  t^{(1+\alpha\beta)(1-\gamma)}\| 
S_{\alpha,\gamma}(t)[u_{0}-h(u(t))]\|
+t^{(1+\alpha\beta)(1-\gamma)}\int_{0}^{t}(t-r)^{
-\alpha\beta-1}[k_1(r)+k_2(r)e^{-\delta r}]dr\\ 
&\leq \sup_{[0,T]}~~ \left(t^{(1+\alpha\beta)(1-\gamma)}\| 
S_{\alpha,\gamma}(t)\|[\Vert u_0\Vert +k]+t^{(1+\alpha\beta)(1-\gamma)}
\int_{0}^{t}(t-r)^{-\alpha\beta-1}[k_1(r)+k_2(r)e^{-\delta r}]dr \right)\\ 
&\leq  r. 
\end{align*}
Hence, $\|\mathcal{P}y\|\leq r$ for any $y\in \Omega_{r}(\mathcal{J})$.
Now, to verify $\mathcal{P}$ is continuous in $\Omega_{r}(\mathcal{J})$, 
let $y_{n}$, $y\in \Omega_{r}(\mathcal{J})$, $n=1,2,\ldots$ 
with $\lim_{n\rightarrow \infty}~y_{n}=y$, that is,
$\lim_{n\rightarrow \infty}~y_{n}(t)=y(t)~$; 
$\text{lim}_{n\rightarrow \infty}~ t^{-(1+\alpha\beta){1-\gamma}}~y_{n}(t)
=t^{-(1+\alpha\beta){1-\gamma}}~y(t)$ and 
$\text{lim}_{n\rightarrow \infty}~t^{-(1+\alpha\beta)(1-\gamma)}y_{n}(t)
=t^{-(1+\alpha\beta)(1-\gamma)}y(t)$ on $\mathcal{J}$. 
Then, $(H1)$ implies that
\begin{align*}
g(t,u_{n}(t),\mathcal{B}(u_{n}(t)))
=&g\left(t,t^{-(1+\alpha\beta)(1-\gamma)}y_{n}(t),
t^{-(1+\alpha\beta)(1-\gamma)}\mathcal{B}(y_{n}(t))\right) \\
&\rightarrow g\left(t,t^{-(1+\alpha\beta)(1-\gamma)}y(t),
t^{-(1+\alpha\beta)(1-\gamma)}\mathcal{B}(y(t))\right)
\end{align*}
as $n\rightarrow \infty$. From $(H2)$, we obtain the inequality 
$$
(t-r)^{-\alpha\beta-1}|g(r,u_{n}(r),\mathcal{B}(u_{n}(r)))|
\leq 2(t-r)^{-(\alpha\beta)(1-\gamma)}[k_1(r)+k_2(r)e^{-\delta r}], 
$$
that is,
\begin{align*} 
\int_{0}^{t}(t-r)^{-\alpha\beta-1}\|g(r,u_{n}(r),\mathcal{B}(u_{n}(r)))
-g(r,u(r),\mathcal{B}(u(r))\|dr~\rightarrow~0
\quad \text{as }~ n\rightarrow \infty.
\end{align*}
Let $t\in [0,T]$. Now,
 \begin{align*}
\|\mathcal{P}y_{n}(t)-\mathcal{P}y(t)\|
\leq t^{(1+\alpha\beta)(1-\gamma)}\Big\|\int_{0}^{t}(t-r)^{\alpha-1}
T_{\alpha}(t-r)(g(r,u_{n}(r),\mathcal{B}(u_{n}(r))-g(r,u(r),\mathcal{B}(u(r)))dr\Big\|.
\end{align*}
Applying Theorem~\ref{Theorem2.1}, we have
\begin{align*}
\|\mathcal{P}y_{n}(t)-\mathcal{P}y(t)\|
\leq C_2t^{(1+\alpha\beta)(1-\gamma)}
\int_{0}^{t}(t-r)^{-\alpha\beta-1}\|g(r,u_{n}(r),
\mathcal{B}(u_{n}(r))-g(r,u(r),\mathcal{B}(u(r))\|dr,
\end{align*}
which tends to 0 as $n\rightarrow \infty$, i.e., 
$\mathcal{P}y_{n}\rightarrow \mathcal{P}y$ pointwise on $\mathcal{J}$. 
Moreover, Theorem~\ref{Theorem3.1} implies that 
$\mathcal{P}y_{n}\rightarrow \mathcal{P}y$ uniformly on $\mathcal{J}$ 
as $n\rightarrow \infty$, that is, $\mathcal{P}$ is continuous.
\end{proof}


\section{Main results}
\label{sec:4}

We prove existence of a mild solution to problem 
\eqref{eq1}--\eqref{eq2} when the associated semigroup 
is compact (Theorem~\ref{Theorem4.1})
and noncompact (Theorem~\ref{Theorem4.2}).

\subsection{Compactness of the semigroup}

Here we assume $Q(t)$ to be compact.

\begin{Theorem}
\label{Theorem4.1} 
Let $-1 <\beta < 0$, $0 < \omega < \frac{\pi}{2}$ and 
$\mathcal{A} \in \Theta_\omega^\beta$. If $ Q(t)(t>0)$ is compact and $(H1)$--$(H4)$ hold, 
then there exists a mild solution of \eqref{eq1}--\eqref{eq2} in $\Omega_r(\mathcal{J})$ 
for every $u_0 \in D(\mathcal{A}^\theta)$ with $\theta > 1 + \beta$.
\end{Theorem}

\begin{proof}
Because we assume $ Q(t)$ to be compact, then the equicontinuity of $ Q(t)(t>0)$ is assured. 
Moreover, by Theorems~\ref{Theorem3.1} and \ref{Theorem3.2}, $\mathcal{P}:\Omega_r(\mathcal{J}) 
\to \Omega_r(\mathcal{J})$ is continuous and bounded and $\varepsilon:\Omega_r(\mathcal{J}) 
\to \Omega_r(\mathcal{J})$ is bounded, continuous, 
and $\{\varepsilon y : y \in \Omega_r(\mathcal{J})\}$ 
is equicontinuous. We can write $\epsilon:\Omega_r(\mathcal{J}) \to \Omega_r(\mathcal{J})$ 
by $(\epsilon y)(t) = (\epsilon^1y)(t)+(\epsilon^2y)(t)$, where
\begin{align*}
(\epsilon^1y)(t)
&=t^{(1+\alpha\beta)(1-\gamma)} S_{\alpha,\gamma}(t)(u_0-h(u)) 
= t^{1+\alpha \beta)(1-\gamma)}\mathcal{I}_{0^+}^{\gamma(1-\alpha)} 
t^{\alpha-1}  T_\alpha(t) (u_0-h(u))\\
&=\frac{t^{(1+\alpha \beta)(1-\gamma)}}{\Gamma(\gamma(1-\alpha))}
\int_0^t (t-r)^{\gamma(1-\alpha)-1}r^{\alpha-1}\int_0^\infty \alpha 
\theta M_\alpha(\theta)   Q(r^\alpha \theta) (u_0-h(u)) d\theta dr\\
&=\frac{\alpha t^{(1+\alpha \beta)(1-\gamma)}}{\Gamma(\gamma(1-\alpha))}
\int_0^t \int_0^\infty (t-r)^{\gamma(1-\alpha)-1}r^{\alpha-1} 
\theta M_\alpha(\theta) Q(r^\alpha \theta) (u_0-h(u)) d\theta dr
\end{align*}
and
\begin{align*}
(\varepsilon^2y)(t)=t^{(1+\alpha\beta)(1-\gamma)} 
\int_0^t (t-r)^{\alpha-1}  
T_\alpha (t-r)g(r,u(r)(\mathcal{B}u)r)dr.
\end{align*}
For $\sigma > 0$ and $ \theta \in (0, t)$, we define an operator 
$\varepsilon_{ \theta,\sigma}^1$ on $\Omega_r(\mathcal{J})$ by
\begin{align*}
(\varepsilon_{ \theta,\sigma}^1y)(t)
&= \frac{t^{(1+\alpha \beta)(1-\gamma)}}{\Gamma(\gamma(1-\alpha))}
\int_\theta^t \int_\sigma^\infty (t-r)^{(1-\alpha)\gamma-1}r^{\alpha-1} 
\theta M_\alpha(\theta)   Q(r^\alpha \theta) (u_0-h(u)) d\theta dr\\
&=\frac{\alpha t^{(1+\alpha \beta)(1-\gamma)}}{\Gamma(\gamma(1-\alpha))} 
\mathfrak{T}( \theta^\alpha \sigma) \int_ \theta^t \int_\sigma^\infty 
(t-r)^{(1-\alpha)\gamma-1}r^{\alpha-1} \theta M_\alpha(\theta) 
Q(r^\alpha \theta- \theta^\alpha \sigma) (u_0-h(u)) d\theta dr.
\end{align*}
Since $\mathfrak{T}(\epsilon^\alpha \delta)$ is compact, 
$\mathcal{V}_{ \theta,\sigma}^1(t)=\{\varepsilon_{ \theta,\sigma}^1 y)(t),y 
\in \Omega_r(\mathcal{J})\}$ is precompact in $\mathcal{Y}$ for all $\theta \in (0,t)$ 
and $\delta > 0$. Moreover, for any $y \in \Omega_r(\mathcal{J})$, one has
\begin{align*}
\|(\varepsilon^1y)(t)-&(\varepsilon_{ \theta,\sigma}^1y)(t)\|\\ 
& \leq \mathcal{K}(\alpha,\gamma) \bigg\|t^{(1+\alpha \beta)(1-\gamma)} 
\int_0^t \int_0^\sigma (t-r)^{\gamma(1-\alpha)-1}r^{\alpha-1} \theta M_\alpha(\theta) 
Q(r^\alpha \theta) (u_0-h(u)) d\theta dr\bigg\| \\
&~~~~+\mathcal{K}(\alpha,\gamma) \bigg\|t^{(1+\alpha \beta)(1-\gamma)} 
\int_0^ \theta \int_\sigma^\infty (t-r)^{\gamma(1-\alpha)-1}r^{\alpha-1} 
\theta M_\alpha(\theta)   Q(r^\alpha \theta) (u_0-h(u)) d\theta dr\bigg\|\\
& \leq \mathcal{K}(\alpha,\gamma) t^{(1+\alpha \beta)(1-\gamma)} 
\int_0^t \int_0^\sigma (t-r)^{\gamma(1-\alpha)-1}r^{\alpha-1} 
\theta M_\alpha(\theta) r^{-\alpha \gamma-\alpha}
\| (u_0-h(u))\| \theta^{-\beta-1} d\theta dr \\
&~~~~+\mathcal{K}(\alpha,\gamma) t^{(1+\alpha \beta)(1-\gamma)} 
\int_0^ \theta \int_\sigma^\infty (t-r)^{\gamma(1-\alpha)-1}r^{\alpha-1} 
\theta M_\alpha(\theta)  r^{-\alpha\beta-\alpha} 
\theta^{-\beta-1} \| (u_0-h(u))\| d\theta dr \\
& = \mathcal{K}(\alpha,\gamma) t^{(1+\alpha \beta)(1-\gamma)} 
\int_0^t (t-r)^{\gamma(1-\alpha)-1}r^{-\alpha\beta-1} \| (u_0-h(u))\| dr 
\int_0^\sigma \theta^{-\beta} M_\alpha(\theta) d\theta \\
&~~~~+  \mathcal{K}(\alpha,\gamma) t^{(1+\alpha \beta)(1-\gamma)} 
\int_0^ \theta (t-r)^{\gamma(1-\alpha)-1}r^{-\alpha \beta-1} \| (u_0-h(u))\| dr 
\int_\eta^\infty \theta^{-\beta} M_\alpha(\theta) d\theta \\
& \leq \mathcal{K} t^{-\alpha\gamma(1+\beta)} \|(u_0-h(u))\| 
\int_0^\eta \theta^{-\beta} M_\alpha(\theta) d\theta\\
&~~~~+  \mathcal{K} t^{-\alpha \gamma(1+\beta)} [\|u_0\|+k] 
\int_0^ \theta (1-s)^{\gamma(1-\alpha)-1}r^{-\alpha \beta-1}dr
\int_\eta^\infty \theta^{-\beta} M_\alpha(\theta) d\theta\\
&\quad \to 0 ~\text{ as } \theta \to 0 \text{ and } ~\sigma \to 0,
\end{align*}
where $\mathcal{K}(\alpha,\gamma)=\frac{\alpha}{\Gamma(\gamma(1-\alpha))}$.
Therefore, $\mathcal{V}_{ \theta,\sigma}^1(t)=\{\varepsilon_{ \theta,\sigma}^1 y)(t),y 
\in \Omega_r(\mathcal{J})\}$ are arbitrarily close to 
$\mathcal{V}^1(t)=\{\varepsilon^1 y)(t),y \in \Omega_r(\mathcal{J})\}$ for $t>0$. 
Hence, $\mathcal{V}^1(t)$, for $t> 0$, is precompact in  $\mathcal{Y}$.
For $ \theta \in (0,t)$ and $\sigma > 0$, we can present an operator 
$\varepsilon_{ \theta,\sigma}^2$ on $\Omega_r(\mathcal{J})$ by
\begin{align*}
(\varepsilon_{ \theta,\sigma}^2y)(t)
=&\alpha t^{(1+\alpha\beta)(1-\gamma)} \int_0^{t- \theta} 
\int_\sigma^\infty \theta M_\alpha(\theta)(t-r)^{\alpha-1} 
Q((t-r)^\alpha \theta)g(r,u(r),(\mathcal{B}u)r)d\theta dr\\
=& \alpha t^{(1+\alpha \beta)(1-\gamma)}\mathfrak{T}( \theta^\alpha \sigma)
\int_0^{t- \theta} \int_\sigma^\infty \theta M_\alpha(\theta)(t-r)^{\alpha-1} 
Q((t-r)^\alpha\theta- \theta^\alpha \sigma)g(r,u(r),(\mathcal{B}u)r)d\theta dr.
\end{align*}
Hence, due to the compactness of $  Q( \theta^\alpha \sigma)$,
$\mathcal{V}_{ \theta,\sigma}^2(t)
=\{\varepsilon_{ \theta,\sigma}^2 y)(t),y \in \Omega_r(\mathcal{J})\}$ 
is precompact in $\mathcal{Y}$ for all  $\theta \in (0,t)$ and $\sigma > 0$. 
For every $y \in \Omega_r(\mathcal{J})$, we get
\begin{align*}
\|(\varepsilon^2y)(t)&-(\varepsilon_{ \theta,\sigma}^2y)(t)\|\\  
\leq& \bigg\|\alpha t^{(1+\alpha\beta)(1-\gamma)}\Big( 
\int_0^t \int_0^\sigma \theta M\alpha(\theta)(t-r)^{\alpha-1}
Q((t-r)^\alpha \theta)g(r,u(r),(\mathcal{B}u)r)d\theta dr\Big)\bigg\|\\
&+\bigg\|\alpha t^{(1+\alpha \beta)(1-\gamma)} \Big(\int_{t- \theta}^t 
\int_\sigma^\infty (t-r)^{\alpha-1}\theta M_\alpha(\theta) 
Q((t-r)^\alpha \theta)g(r,u(r),(\mathcal{B}u)r)d\theta dr\Big)\bigg\|\\
\leq & \alpha  C_0 t^{(1+\alpha \beta)(1-\gamma)} \Big(\int_0^t 
(t-r)^{-\alpha \beta - 1}[k_1(r)+k_2(r)e^{-\delta r}] dr 
\int_0^\sigma\theta^{-\beta} M_\alpha(\theta) d\theta\Big)\\
&~~~~+ \alpha C_0 t^{(1+\alpha \beta)(1-\gamma)}\Big(\int_{t
- \theta}^t (t-r)^{-\alpha\beta - 1}[k_1(r)+k_2(r)e^{-\delta r}] dr 
\int_0^\infty \theta^{-\beta} M_\alpha(\theta) d\theta\Big)\\
\leq & \alpha  C_0 t^{(1+\alpha\beta)(1-\gamma)} \Big(\int_0^t 
(t-r)^{-\alpha \beta - 1}[k_1(r)+k_2(r)e^{-\delta r}] dr 
\int_0^\sigma \theta^{-\beta} M_\alpha(\theta) d\theta\Big)\\
&~~~~+\frac{\alpha  C_0 \Gamma(1-\beta)}{\Gamma(1-\alpha \beta)} 
t^{(1+ \alpha \beta)(1-\gamma)} \Big(\int_{t- \theta}^t 
(t-r)^{-\alpha \beta - 1}[k_1(r)+k_2(r)e^{-\delta r}] dr\Big)
\to 0~\text{as}~\sigma \to 0.
\end{align*}
Therefore, $\mathcal{V}_{ \theta,\sigma}^2(t)
=\{\varepsilon_{ \theta,\sigma}^2 y)(t),y \in \Omega_r(\mathcal{J})\}$ 
are arbitrarily close to $\mathcal{V}^2(t)=\{\varepsilon^2 y)(t),y 
\in \Omega_r(\mathcal{J})\}$, $t > 0$. This implies the relative compactness 
of $\mathcal{V}^2(t)$, $t > 0$, in $\mathcal{Y}$. Also, 
$\mathcal{V}(t)=\{\varepsilon y)(t),y \in \Omega_r(\mathcal{J})\}$ 
is relatively compact in $\mathcal{Y}~\forall t \in [0,T]$. It follows,
from the Arzela--Ascoli theorem, that
$\{\varepsilon y,y \in \Omega_r(\mathcal{J})\}$ is relatively compact
for $\varepsilon$, it is continuous, and $\{\varepsilon y,y \in \Omega_r(\mathcal{J})\}$ 
is relatively compact. This implies, by the Schauder fixed point theorem, existence of 
a fixed point $y^* \in \Omega_r(\mathcal{J})$ of $\varepsilon$. 
Let $u^*(t) := t^{(1+\alpha \beta)(\gamma - 1)}y^*(t)$. Then, 
$u^*$ is a mild solution of \eqref{eq1}--\eqref{eq2}.
\end{proof} 
 

\subsection{Non-compactness of the semigroup}

Now, we assume that $Q(t)$ is noncompact. We need the following supplementary condition:
\begin{enumerate}
\item[$(H5)$] There exists a constant $k > 0$ satisfying
\begin{align*}
\Theta(g(t,\mathfrak{E}_{1},\mathfrak{E}_{2})) 
\leq k \Theta(\mathfrak{E}_{1},\mathfrak{E}_{2})~\text{for a.a. } t \in [0,T]
\end{align*}
and for every bounded subsets $\mathfrak{E}_{1},\mathfrak{E}_{2} \subset  \mathcal{Y}$.
\end{enumerate}
 
\begin{Theorem}
\label{Theorem4.2} 
Let $-1 < \beta < 0$, $0 < \omega < \frac{\pi}{2}$, and $\mathcal{A} \in \Theta_\omega^\beta$. 
Suppose $(H1)$--$(H5)$ hold. Then, \eqref{eq1}--\eqref{eq2} has a mild solution in 
$\Omega_r(\mathcal{J})$ for every $u_0 \in D(g^\theta)$ with $\theta > 1 + \beta$.
\end{Theorem}

\begin{proof}
From Theorems~\ref{Theorem3.1} and \ref{Theorem3.2}, we get 
that $\varepsilon:\Omega_r(\mathcal{J}) \to \Omega_r(\mathcal{J})$ is continuous, bounded, 
and $\{\varepsilon y: y \in \Omega_r(\mathcal{J})\}$ is equicontinuous. Also, we prove 
that there is a subset of $\Omega_r(\mathcal{J})$ such that $\varepsilon$ is compact in it.
For any bounded set $\mathbb{P}_0 \subset \Omega_r(\mathcal{J})$, set
\begin{align*}
\varepsilon^{(1)}(\mathbb{P}_0) = \varepsilon(\mathbb{P}_0),\varepsilon^{(n)}(\mathbb{P}_0) 
= \varepsilon(\bar{co}(\varepsilon^{(n-1)}(\mathbb{P}_0))), 
\quad n = 2,3,\ldots
\end{align*}
For any $\epsilon > 0$, we can obtain from Propositions~\ref{Proposition2.2}--\ref{Proposition2.4} 
a subsequence $\{y_n^{(1)}\}_{n=1}^\infty \subset \mathbb{P}_0$ satisfying
\begin{align*}
\Theta(\varepsilon^{(1)}(\mathbb{P}_0(t))) 
&\leq 2 \Theta \bigg(t^{(1+\alpha \beta)(1-\gamma)} \int_0^t (t-r)^{\alpha-1}  
T_\alpha(t-r) g(r,\{r^{-(1+\alpha \beta)(1-\gamma)} (y_n^{(1)}(r),
\mathcal{B}y_n^{(1)}(r))\}_{n=1}^\infty )dr \bigg)\\
&\leq 4 C_p t^{(1+\alpha \beta)(1-\gamma)}\Big(\int_0^t (t-r)^{-\alpha \beta-1} 
\Theta ( g(r,\{r^{-(1+\alpha \beta)(1-\gamma)} (y_n^{(1)}(r),
\mathcal{B}y_n^{(1)}(r))\}_{n=1}^\infty ))dr \Big)\\
&\leq  4 C_p kt^{(1+\alpha \beta)(1-\gamma)}\Theta (\mathfrak{P}_0)\Big(
\int_0^t (t-r)^{-\alpha \beta-1} r^{-(1+\alpha\beta)(1-\gamma)} dr \Big)\\
&= 4 C_p kt^{-\alpha \beta}\Theta (\mathfrak{P}_0) \Big(\frac{\Gamma(
-\alpha \beta) \Gamma((-\alpha \beta 
+ \gamma(1+\alpha \beta))}{\Gamma(-2\alpha \beta + \gamma(1+\alpha \beta))}\Big).
\end{align*}
Since $\epsilon$ is arbitrary, then 
\begin{align*}
\Theta(\varepsilon^{(1)}(\mathbb{P}_0(t))) 
\leq 4 C_p kt^{-\alpha \beta}\Theta (\mathbb{P}_0) \Big(\frac{\Gamma (-\alpha \beta) 
\Gamma((-\alpha \beta + \gamma(1+\alpha \beta))}{\Gamma(-2\alpha\beta 
+ \gamma(1+\alpha \beta))}\Big).
\end{align*}
Again, for any $\epsilon > 0$, we can get from Propositions~\ref{Proposition2.2}--\ref{Proposition2.4} 
a subsequence $\{y_n^{(2)},\mathcal{B}y_{n}^{(2)}\}_{n=1}^\infty 
\subset \bar{co}(\varepsilon^{(1)}(\mathbb{P}_0))$, which implies that
\begin{align*}
\Theta(\varepsilon^{(2)}(\mathbb{P}_0(t))) 
&= \Theta(\varepsilon (\bar{co}(\varepsilon^{(1)}(\mathbb{P}_0(t))))) \\
&\leq 2 \Theta \bigg(t^{(1+\alpha\beta)(1-\gamma)} 
\int_0^t (t-r)^{\alpha-1} \mathcal{Q}_\alpha(t-r)  
g(r,\{r^{-(1+\alpha\beta)(1-\gamma)} (y_n^{(2)}(r),
\mathcal{B}y_n^{(2)}(r))\}_{n=1}^\infty )dr\bigg)\\
& \leq 4 C_p t^{(1+\alpha\beta)(1-\gamma)}\Big(
\int_0^t (t-r)^{-\alpha \beta-1} \Theta ( g(r,\{r^{-(1+\alpha \beta)(1-\gamma)}
(y_n^{(2)}(r),\mathcal{B}y_n^{(2)}(r))\}_{n=1}^\infty )dr \Big)\\
& \leq 4 C_p kt^{(1+\alpha\beta)(1-\gamma)}\Big(\int_0^t 
(t-r)^{-\alpha\beta-1} \Theta (r^{-(1+\alpha \beta)(1-\gamma)}(
\{y_n^{(2)}(r),\mathcal{B}y_{n}^{(2)}(r)\}_{n=1}^\infty))dr\Big)\\
& \leq 4 C_p kt^{(1+\alpha\beta)(1-\gamma)}\Big(\int_0^t 
(t-r)^{-\alpha \beta-1} r^{-(1+\alpha \beta)(1-\gamma)} \Theta(\{
y_n^{(2)}(r),\mathcal{B}y_{n}^{(2)}\}_{n=1}^\infty)dr\Big)\\
& \leq \frac{(4 C_p k)^2 t^{(1+\alpha \beta)(1-\gamma)} 
\Gamma(-\alpha \beta)\Gamma(-\alpha\beta + \gamma(1+\alpha \beta))}{\Gamma(-2 
\alpha \beta + \gamma(1+\alpha \beta))} \Theta (\mathbb{P}_0)\\
&~~~~\times\Big( \int_0^t (t-r)^{-\alpha \beta-1} r^{-(1+\alpha \beta)(1-\gamma)
-\alpha \beta} dr \Big)\\
&= \Big(\frac{(4 C_p k)^2 t^{-2\alpha \beta} \Gamma^2(-\alpha \beta)
\Gamma(-\alpha\beta + \gamma(1+\alpha \beta))}{\Gamma(-3\alpha \beta 
+ \gamma(1+\alpha \beta))}\Big)\Theta (\mathbb{P}_0).
\end{align*}
Now,
\begin{align*}
\Theta(\varepsilon^{(n)}(\mathbb{P}_0(t))) 
\leq \frac{(4 C_p k)^n t^{-n\alpha \beta} 
\Gamma^n(-\alpha \beta)\Gamma(-\alpha \beta 
+ \gamma(1+\alpha \beta))}{\Gamma(-(n+1) \alpha\beta 
+ \gamma(1+\alpha\beta))}\Theta (\mathbb{P}_0),
\quad n \in \mathbb{N}.
\end{align*}
Let $M = 4 C_pkT^{-\alpha \beta} \Gamma(-\alpha \beta)$. 
We can find $m,k \in \mathbb{N}$ big enough such that 
$\frac{1}{k}<\alpha \beta< \frac{1}{k-1}$ and $\frac{n+1}{k}>2$ 
for $n > m \Gamma(-(n+1)\alpha \beta + \gamma(1+\alpha \beta))
>\Gamma(\frac{n+1}{k})$, that is,
\begin{align*}
\frac{(4 C_pk)^n T^{-n \alpha \beta} \Gamma^n(-\alpha \beta)
\Gamma(-\alpha \beta + \gamma(1+\alpha \beta))}{\Gamma(-(n+1)\alpha \beta 
+ \gamma(1+\alpha\beta))} < \frac{(4 C_pk)^n T^{-n \alpha \beta} \Gamma^n(
-\alpha \beta)\Gamma(-\alpha \beta + \gamma(1
+\alpha \beta))}{\Gamma\bigg(\frac{n+1}{k}\bigg)}.
\end{align*}
Replacing $(n+1)$ by $(j+1)k$, then the right-hand side of the
inequality given above becomes
\begin{align*}
\frac{M^{(j+1)k-1}\Gamma(-\alpha\beta + \gamma(1+\alpha \beta))}{\Gamma(j+1)} 
= \frac{(M^k)^j M^{k-1}\Gamma(-\alpha \beta + \gamma(1+\alpha \beta))}{j!} 
\to 0~\text{ as }~j \to \infty.
\end{align*}
Therefore, there exists a constant $n_0 \in \mathbb{N}$ such that
\begin{align*}
\frac{(4 C_pk^n t^{-n \alpha \beta} \Gamma^n(-\alpha \beta)
\Gamma(-\alpha \beta + \gamma(1+\alpha \beta)}{\Gamma(-(n+1) 
\alpha \beta + \gamma(1+\alpha \beta))} 
\leq \frac{(4 C_pk)^{n_0} T^{-n_0 \alpha \beta} 
\Gamma^{n_0}(-\alpha \beta)\Gamma(-\alpha \beta 
+ \gamma(1+\alpha \beta))}{\Gamma(-(n_0^+1)\alpha\beta 
+ \gamma(1+\alpha \beta))} = p <1.
\end{align*}
Now, $\Theta(\varepsilon^{(n_0)}(\mathbb{P}_0(t))) \leq p \Theta(\mathbb{P}_0)$.
Since $\varepsilon^{(n_0)}(\mathbb{P}_0(t))$ is bounded and equicontinuous, 
it follows from Proposition~\ref{Proposition2.2} that
\begin{align*}
\Theta(\varepsilon^{(n_0)}(\mathbb{P}_0)) 
= \max_{t \in [0,T]} \Theta (\varepsilon^{n_0}(\mathbb{P}_0(t))).
\end{align*}
Hence, $\Theta(\varepsilon^{n_0}(\mathbb{P}_0)) \leq p \Theta(\mathbb{P}_0)$,
where $p < 1$. Using a similar technique as in Theorem~\ref{Theorem4.1}, 
we obtain $C$ in $\Omega_r(\mathcal{J})$ with $\varepsilon( C) \subset  C$ and 
$\varepsilon( C)$ compact. By applying Schauder's fixed point theorem, 
we  obtain a fixed point $y^*$ in $\Omega_r(\mathcal{J})$ of $\varepsilon$. 
Let $u^*(t) = t^{(1+\alpha \beta)(\gamma - 1)} y^*(t)$. Then, $u^*(t)$  
is a mild solution of \eqref{eq1}--\eqref{eq2}.
\end{proof} 
 

\section{Example}
\label{sec:5}

As an illustrative example, let us consider the following Hilfer fractional 
partial differential equation with a nonlocal condition:
\begin{equation}
\label{eq4}
\begin{gathered}
D^{\frac{3}{4},\frac{1}{2}}_{0^+}z(t, y)+\partial^{2}_{x}z(t, y)
=G(t, z(t, y), \mathcal{B}z(t, y)),
\quad t\in (0, 1], \quad y\in [0, \pi],\\
z(t, 0)=z(t, \pi)=0,
\quad t\in (0, 1],\\
I^{\frac{1}{8}}_{0^+}[z(t, y)]\vert_{t=0}
+\sum\limits_{i=1}^{m}c_{i}z(t_{i}, y)
=z_{0}(y), \quad y\in [0, \pi],
\end{gathered}
\end{equation}
where $\alpha=\frac{3}{4}$, $\gamma=\frac{1}{2}$, 
$0< t_{1}<\cdots<t_{m}< 1$, and $c_{i}$, $i=1,\ldots, m$,
are given constants. Let us take the nonlinear function $G(t, z(\cdot), \mathcal{B}z(\cdot))
=y\cos z(t, y)+\int_0^te^{-(t-s)}\sin z(t, y)ds$ and the nonlocal function 
$h(z(t, \cdot))=\sum_{i=1}^{m}c_{i}z( t_{k}, \cdot)$. Assume that $\mathcal{Y}=L^{2}[0, \pi]$ 
and define $\mathcal{A}: D(\mathcal{A})\subset \mathcal{Y}\rightarrow \mathcal{Y}$ 
by $\mathcal{A}z=z_{yy}$ with domain
$$
D(\mathcal{A})=\left\lbrace z\in \mathcal{Y}: ~z_{y}, z_{yy}
\in \mathcal{Y}, ~ z(t, 0)=z(t, \pi)=0\right\rbrace.
$$
It follows from the work in \cite{p1} that there exists constants $\delta,\epsilon>0$ such that 
$\mathcal{A}+\delta\in \odot_{\frac{\pi}{2}-\epsilon}^{\frac{\pi}{2}-1}(\mathcal{Y})$.
It is known that $\mathcal{A}$ is the infinitesimal generator of a differentiable semigroup 
$Q(t)(t>0)$ in $\mathcal{Y}$ given by
\begin{equation*}
(Q(t)x)(y)=
\begin{cases}
\int_{0}^{\pi}\psi(t, y-s)x(s)ds, ~t>0,\\
x(y), \quad t=0,
\end{cases}
\end{equation*}
where 
$$
\psi(t, y)=\frac{1}{\sqrt{4\pi t}}e^{-\frac{y^{2}}{4t}}, 
\quad t>0, \quad 0<y<\pi,
$$
and $x(t)(y)=z(t, y)$. This implies $\Vert Q(t)\Vert\leq1$ and leads to its compactness property. 
We can check that all hypotheses $(H1)$--$(H4)$ are fulfilled. Hence, our Theorem~\ref{Theorem4.1} 
can be applied ensuring that problem \eqref{eq4} admits a mild solution.


\section{Conclusion} 
\label{sec:6}

In this paper, we applied Schauder's fixed point theorem to investigate the solvability 
of a class of Hilfer fractional integro-differential equations involving almost sectorial operators. 
We discussed both cases of compactness and noncompactness, related to associated semigroup operators. 
The obtained existence results were subject to an appropriate set of sufficient conditions. 
As a future direction of research, it is desirable to consider the study of $\psi$-Hilfer fractional 
nonlocal nonlinear stochastic systems involving almost sectorial operators and impulsive 
effects, generalizing the current work. Another open line of research 
consists to develop numerical methods to approximate the mild solutions 
predicted by our Theorems~\ref{Theorem4.1} and \ref{Theorem4.2}.


\authorcontributions{Conceptualization, A.~Debbouche; 
methodology, D.~F.~M. Torres; validation, A.~Debbouche; 
formal analysis, D.~F.~M. Torres; investigation, D.~F.~M.~Torres; 
writing--original draft preparation, K.~Karthikeyan; 
writing--review and editing, A.~Debbouche and D.~F.~M.~Torres; 
supervision, A.~Debbouche; project administration, D.~F.~M.~Torres.}


\funding{Debbouche and Torres were supported by FCT within project UIDB/04106/2020 (CIDMA).}


{\conflictsofinterest{The authors declare no conflict of interest. 
The funder had no role in the design of the study; in the collection, 
analyses, or interpretation of data; in the writing of the manuscript, 
or in the decision to publish the results.}


\reftitle{References}


\end{document}